\documentclass[11pt,twoside]{article}

\date{}
\usepackage{epsfig}

\begin{document}

%\begin{figure}[h]
%\includegraphics[scale=1.40]{ARM.png}
%{\small{\bf Vol. x, No. x (200x) xxx-xxx}}
%\end{figure}

\begin{center}
{\Large \bf Affine Geometry of Space Curves}

{\small \bf Mehdi Nadjafikhah$^{a,}$\footnote{\footnotesize
Corresponding Author. E-mail Address: m\_nadjafikhah@iust.ac.ir},
Ali Mahdipour--Shirayeh$^b$}
\end{center}
{\tiny \bf $^{a,b}$Department of Mathematics, Iran University of Science and Technology, Narmak, Tehran 16846-13114, IRAN}\\

\newtheorem{Theorem}{\quad Theorem}[section]
\newtheorem{Definition}[Theorem]{\quad Definition}
\newtheorem{Corollary}[Theorem]{\quad Corollary}
\newtheorem{Lemma}[Theorem]{\quad Lemma}
\newtheorem{Example}[Theorem]{\quad Example}

\begin{abstract}
This paper is devoted to the complete classification of space
curves under affine transformations in the view of Cartan's
theorem. Spivak has introduced the method but has not found the
invariants. Furthermore, for the first time, we propound a
necessary and sufficient condition for the invariants. Then, we
study the shapes of space curves with constant curvatures in
detail and suggest their applications in physics, computer vision
and image processing.
\\
\copyright {\small \hspace{0.15cm}200x Published by Islamic AZAD
University-Karaj Branch.}
\end{abstract}

\section{Introduction}
Classification of curves has a significant place in geometry,
physics, mechanics, computer vision and image processing. In
geometrical sense, a plane curve with constant curvature, up to
special affine transformations may be either an ellipse, a
parabola or a hyperbola \cite{Sp}. This classification will be
obtained by the concept of invariants. Geometry of curves in
spaces with dimension $\geq 3$ has studied with geometers such as
Guggenheimer \cite{Gu}, Spivak \cite{Sp} and etc. The aim was
finding the invariants of curves under transformations. On the
other hand, in \cite{Sp}, study of space curves in the view of
Cartan's theorem was started but has not completed yet.

This paper can be viewed as a continuation of the work \cite{Sp},
where the authors began the classification of space curves up to
special affine transformations. We determine all of differential
invariants and our method is different from the method of
Guggenheimer and other existing methods. Also, for the first time,
we prove a necessary and sufficient condition for the invariants
in order that complete the classification. Moreover, we classify
the shapes of space curves of constant curvatures which has a wide
variety of applications in physics, computer vision and image
processing. The general form of these shapes are exist in
\cite{Gu}, but here we try to discuss them in more details.

In physics, classification of curves up to affine transformations
has a special position in the study of rigid motions. Suppose we
have a particle moving in 3D space and that we want to describe
the trajectory of this particle. Especially, each curve in a three
dimensional space could be imagined as a trajectory of a particle
with a specified mass in the view of an observer. By
classification of curves we can, in fact, obtain conservation
laws.

Computer vision deals with image understanding at various levels.
At the low level, it addresses issues such us planar shape
recognition and analysis. Some results on differential invariants
associated to space curves are relevant to space object
recognition under different views and partial occlusion. The
evolution of space shapes under curvature controlled diffusion
have applications in geometric shape decomposition, smoothing, and
analysis, as well as in other image processing applications (see,
e.g. \cite{MP,MWTN}) and similar to recent results for planer
shapes. For instance, there are some important applications of
moving frames method in use of the differential invariant
signatures \cite{Ol1}. In \cite{Ca-Ol,Ca-Ol1} there exist some
applications to the problem of object recognition and symmetry.
Also, joint differential invariants has been proposed as
noise-resistant alternatives to differential invariant signatures
in computer vision \cite{Ca-Mo}. Practical applications of the
derived shapes in the latest section are related to invariant
signatures, object recognitions, and symmetry of 3D shapes via the
generalization of them from 2D shapes to 3D ones.

In the next section, we state some preliminaries about
Maurer-Cartan forms and a way of classification of maps with the
notable role of Maurer-Cartan forms and Cartan's theorem. In
section three, classification of space curves in ${\bf R}^3$ under
the action of affine transformations is discussed. Finally, in the
last section, we study the shapes of space curves with constant
curvatures and propose some applications of these shapes in
physics, computer vision and image processing.
%
%
%%%%%%%%%%%%%%%%%%%%%%%%%%%%%%%%%%%%%%%%%%%%%%%%%%%%%%%%%%%%%%%%%
\section{Maurer-Cartan form}
Let $G\subset{\rm GL}(n,{\bf R})$ be a matrix Lie group with Lie
algebra ${\mathcal L}$ and $P:G\rightarrow{\rm Mat}(n\times n)$ be
a matrix-valued function which embeds $G$ into ${\rm Mat}(n\times
n)$, the vector space of $n\times n$ matrices with real entries.
Its differential is $dP_B: T_BG\rightarrow T_{P(B)}{\rm
Mat}(n\times n)\simeq {\rm Mat}(n\times n)$.
\paragraph{Definition 2.1} The 1-form $\omega_B=\{P(B)\}^{-1}\cdot dP_B$
of $G$ is called the {\it Maurer-Cartan form}. It is often written
$\omega = P^{-1}\cdot dP$. The Maurer-Cartan form is in fact the
unique left invariant ${\mathcal L}-$valued 1-form on $G$ such
that $\omega_{\rm Id}:T_{\rm Id}G\rightarrow{\mathcal L}$ is the
identity map. The Maurer-Cartan form $\omega$ satisfies in {\it
Maurer-Cartan equation} $d\,\omega = - \omega\,\wedge\,\omega.$
The Maurer-Cartan form is the key to classifying maps into
homogeneous spaces of $G$. This process needs to the following
theorem (for a proof we refer the reader to \cite{Iv-La}):
\paragraph{Theorem 2.2 (Cartan)}{\it Let $G$ be a matrix Lie group with Lie
algebra ${\mathcal L}$ and Maurer-Cartan form $\omega$. Let $M$ be
a manifold on which there exists a ${\mathcal L}-$valued 1-form
$\phi$ satisfying $d\,\phi = -\phi\wedge\phi$. Then for any point
$x\in M$ there exist a neighborhood $U$ of $x$ and a map
$f:U\rightarrow G$ such that $f^{\ast}\,\omega = \phi$. Moreover,
Any two such maps $f_1, f_2$ must satisfy $f_1 = L_B\circ f_2$ for
some fixed $B\in G$ ($L_B$ is the left action of $B$ on $G$).}
\paragraph{Corollary 2.3}{\it Given maps $f_1, f_2:M\rightarrow
G$, then $f_1^{\ast}\,\omega =f_2^{\ast}\,\omega$, that is, this
pull-back is invariant, if and only if $f_1 = L_B\circ f_2$ for
some fixed $B\in G$.}\\

In Theorem 2.2, if $M$ is connected and simply-connected, then the
desired map $f$ may be extended to all of $M$ \cite{Wa}. We
suppose that $G$ be the special linear group ${\rm SL}(3,{\bf R})$
as a Lie group and we denote its Lie algebra with ${\rm sl}(3,{\bf
R})$. This Lie group is not simply-connected, so our achievements
are local.
\paragraph{Definition 2.4} An {\it affine transformation} of the Euclidean
space ${\bf R}^3$ is the composition of a translation in ${\bf
R}^3$ among with an element of the general linear group ${\rm
GL}(3,{\bf R})$. An affine transformation is called {\it special}
or {\it unimodular}, if its matrix part is an element of ${\rm
SL}(3,{\bf R})$. The group of special affine transformations is
the connected coefficient, closed subgroup of the Lie group of
affine transformations.\\

The following section is devoted to the study of the properties of
space curves' invariants under the action of volume--preserving
affine transformations, i.e., the special affine group. The number
of essential parameters (dimension of the Lie algebra) is 11. The
natural assumption of differentiability is ${\mathcal C}^5$.

In the next section, by defining the new curve $\alpha_c$ instead
of a considered regular smooth curve $c$, we will see that the
classification of curves in ${\bf R}^3$ and in the viewpoint of
Theorem 3.1 is equivalent to the ones in ${\rm SL}(3,{\bf R})$.
Thus we find the Maurer-Cartan form of ${\rm SL}(3,{\bf R})$ and
then its pull-back via the matrix--valued curve $\alpha_c$. In
fact, $\alpha_c$~s play the role of $f_i$~s in Corollary 2.3. This
tends to a complete set of invariants of $\alpha_c$ as 1-forms on
${\bf R}$. The derived invariants in a corresponding manner
determines curves of ${\rm SL}(3,{\bf R})$. Finally in Theorem
3.3, we find that these invariants also provide a necessary and
sufficient condition for specifying curves in ${\bf R}^3$ when we
supposed the action of special affine group ${\rm SL}(3,{\bf R})$
on ${\bf R}^3$.
%
%%%%%%%%%%%%%%%%%%%%%%%%%%%%%%%%%%%%%%%%%%%%%%%%%%%%%%%%%%%%%%%
\section{Classification of space curves}
In the present section, we achieve the invariants of a space curve
up to special affine transformations. From Theorem 2.2, two curves
in ${\bf R}^3$ are equivalent under special affine
transformations, if they differ with a left action introduced by
an element of ${\rm SL}(3,{\bf R})$ and then a translation.

Let $c:[a,b]\rightarrow{\bf R}^3$ be a curve in three dimensional
space which we call the {\it space curve}, be of class ${\cal
C}^5$ and
\begin{eqnarray}
\det(c',c'',c''')\neq 0,\label{(1)}
\end{eqnarray}
for any point of the domain, that is, we assume that $c'$, $c''$
and $c'''$ are linear independent. Otherwise, if for example,
$c'''$ depends to $c'$ and $c''$ for some interval $[a,b]$, then
we can simply observe that the curve $c$ will sit in ${\bf R}^2$,
which is not our main topic of investigation. Moreover, we can
assume that $\det(c',c'',c''')>0$ for being avoid writing the
absolute value in calculations.\par
For the curve $c$, we consider a new curve, namely
$\alpha_c(t):[a,b]\rightarrow {\rm SL}(3,{\bf R})$, defined by
\begin{eqnarray*}
\alpha_c(t):=\frac{(c',c'',c''')}{\{\det(c',c'',c''')\}^{1/3}}
\end{eqnarray*}
which is well defined on the domain of $c$ into the special linear
group ${\rm SL}(3,{\bf R})$. We can study the new curve in respect
to special affine transformations, i.e. the action of special
affine transformations on first, second and third differentiations
of $c$. If we assume that $A$ is a three dimensional special
affine transformation, then we have the unique representation
$A=\tau\circ B$ which $B$ is an element of ${\rm SL}(3,{\bf R})$
and $\tau$ is a translation in ${\bf R}^3$. If two curves $c$ and
$\bar{c}$ be the same up to an $A$, $\bar{c}=A\circ c$, then we
have
\begin{eqnarray*}
\bar{c}'=B\circ c', \hspace{1cm} \bar{c}''=B\circ c'',
\hspace{1cm} \bar{c}'''=B\circ c'''.
\end{eqnarray*}
Also from $\det B=1$ we obtain
\begin{eqnarray*}
\det(\bar{c}',\bar{c}'',\bar{c}''')&=&\det(B\circ c',B\circ
c'',B\circ c''')=\det(B\circ(c',c'',c'''))\\
&=&\det(c',c'',c''').
\end{eqnarray*}
and so we conclude that $\alpha_{\bar{c}}(t)=B\circ\alpha_c(t)$
and $\alpha_{\bar{c}}=L_B\circ\alpha_c$, where $L_B$ is the left
translation for $B\in{\rm SL}(3,{\bf R})$.

This condition is also necessary because when $c$ and $\bar{c}$
are two space curves in which $\alpha_{\bar{c}}=L_B\circ\alpha_c$
for an element $B\in {\rm SL}(3,{\bf R})$, then we can write
\begin{eqnarray*}
\alpha_{\bar{c}} &  = & \{\det(\bar{c}',\bar{c}'',\bar{c}''')\}^{-1/3} \,(\bar{c}',\bar{c}'',\bar{c}''') \\
& = & \{\det(B\circ(c',c'',c'''))\}^{-1/3}   \, B\circ(c',c'',c''') \\
& = & \{\det(c',c'',c''')\}^{-1/3}\,B\circ(c',c'',c''').
\end{eqnarray*}
Thus $\bar{c}'=B\circ c'$ and there is a translation $\tau$ such
that $A=\tau\circ B$ and $\bar{c}=A\circ c$ where $A$ is a three
dimensional affine transformation. Therefore we have
\paragraph{Theorem 3.1.} {\it Let $c$ and $\bar{c}$ are two space curves. $c$ and $\bar{c}$ are the same
with respect to special affine transformations, i.e.
$\bar{c}=A\circ c$ when $A=\tau\circ B$ for translation $\tau$ in
${\bf R}^3$ and $B\in{\rm SL}(3,{\bf R})$ if and only if
$\alpha_{\bar{c}}=L_B\circ\alpha_c$ where $L_B$ is a left
translation generated by $B$.}\\

From Cartan's theorem, a necessary and sufficient condition for
$\alpha_{\bar{c}}=L_B\circ\alpha_c$ ($B\in {\rm SL}(3,{\bf R})$)
is that for any left invariant 1-form $\omega^i$ on ${\rm
SL}(3,{\bf R})$ we have
$\alpha_{\bar{c}}^{\ast}(\omega^i)=\alpha_{c}^{\ast}(\omega^i)$.
It is equivalent to
$\alpha_{\bar{c}}^{\ast}(\omega)=\alpha_{c}^{\ast}(\omega)$ for
natural ${\rm sl}(3,{\bf R})$-valued 1-form $\omega=P^{-1}\,.\,dP$
for matrix-valued function $P$ which embeds ${\rm SL}(3,{\bf R})$
into ${\rm Mat}{3\times 3}$, the vector space of $3\times 3$
matrices with real entries, and $\omega$ is the Maurer-Cartan
form. \par
We must compute $\alpha_c^{\ast}(P^{-1}.dP)$ which is invariant
under special affine transformations. Its entries are in fact
invariant functions of space curves. It is a $3\times3$ matrix
form which arrays are multiplications of $dt$ (1-forms on
$[a,b]$).

Since $\alpha_c^{\ast}(P^{-1}\cdot dP)=\alpha_c^{-1}\cdot
d\alpha_c$, so we calculate the matrix
$\alpha_c^{-1}\cdot\alpha'_c$  and then multiply it by $dt$ to
have $\alpha_c^{\ast}(P^{-1}\cdot dP)$. we have
\begin{eqnarray*}
\alpha_c^{-1} \!\!\!\!&=&\!\!\!\! \det(c',c'',c''')^{1/3}\!\cdot(c',c'',c''')^{-1} \\
&=&\!\!\! \det(c',c'',c''')^{-2/3}\!\cdot\left(\!\!
\begin{array}{ccc}
 c''_2c'''_3-c''_3c'''_2 & c''_3c'''_1-c''_1c'''_3 &
 c''_1c'''_2-c''_2c'''_1\\c'_3c'''_2-c'_2c'''_3
 & c'_3c'''_5-c'_1c'''_3 & c'_2c'''_1-c'_1c'''_2 \\
c'_2c''_3-c'_3c''_2
 &c'_3c''_1
c'_1c''_3 & c'_1c''_2-c'_2c''_1
\end{array}\!\!\right)
\end{eqnarray*}
which $c=(c_1,c_2,g_3)^T$ as a column matrix be the vector
representation of curve $c$. We also have
\begin{eqnarray*}
[\det(c',c'',c''')]' &=& \det(c'',c'',c''')
+\det(c',c''',c''') + \det(c',c'',c'''') \\
&=&\det(c',c'',c'''').
\end{eqnarray*}
Thus we see that
\begin{eqnarray*}
\alpha'_c=\det(c',c'',c''')^{-1/3}\!\cdot\!\left(\!\!
\begin{array}{ccc}
c''_1 & c'''_1  & c''''_1 \\
c''_ 2& c'''_ 2 & c''''_ 2\\
c''_3 & c'''_ 3 & c''''_3
\end{array}\!\!\right)
\!\!- \frac{1}{3}\det(c',c'',c''')^{-4/3}\!\cdot\!
\left(\!\!\begin{array}{ccc}
c'_1 & c''_1  & c'''_1 \\
c'_2 & c''_2 & c'''_2\\
c'_3 & c''_3 & c'''_3
\end{array}\!\!\right)\!.
\end{eqnarray*}
After some computations, finally we find that
$\alpha_c^{-1}\cdot\alpha'_c$ is in the following multiple of $dt$
\begin{eqnarray*}\left(
\begin{array}{ccc}
-\frac{\det(c',c'',c'''')}{3\det(c',c'',c''')}& 0
& \frac{\det(c'',c''',c'''')}{\det(c',c'',c''')}\\[0.5mm]
1 & -\frac{\det(c',c'',c'''')}{3\det(c',c'',c''')}
& -\frac{\det8c',c''',c'''')}{\det(c',c'',c''')}\\[0.5mm]
0 & 1 & \frac{2\det(c',c'',c'''')}{3\det(c',c'',c''')}
\end{array}\right).
\end{eqnarray*}
Clearly, the trace of the last matrix is zero and entries of
$\alpha_c^{\ast}(P^{-1}\cdot dP)$ and therefore entries of the
above matrix, are invariants of the group action.\par
Therefore according to Theorem 3.1, two space curves $c,
\bar{c}:[a,b]\rightarrow{\bf R}^3$ are the same under special
affine transformations if we have
\begin{eqnarray*}
\frac{\det(c',c'',c'''')}{\det(c',c'',c''')}&=&
\frac{\det(\bar{c}',\bar{c}'',\bar{c}'''')}{\det(\bar{c}',\bar{c}'',\bar{c}''')}\\[1mm]
\frac{\det(c'',c''',c'''')}{\det(c',c'',c''')} & =&
\frac{\det(\bar{c}'',\bar{c}''',\bar{c}'''')}{\det(\bar{c}',\bar{c}'',\bar{c}''')}\\[1mm]
\frac{\det(c',c''',c'''')}{\det(c',c'',c''')}&=&
\frac{\det(\bar{c}',\bar{c}''',\bar{c}'''')}{\det(\bar{c}',\bar{c}'',\bar{c}''')}.
\end{eqnarray*}
\medskip We may use of a proper parametrization
$\sigma:[a,b]\rightarrow [0,l]$ such that the parameterized  curve
$\gamma=c\circ\sigma^{-1}$ satisfies in condition
$\det(\gamma'(s),\gamma''(s),\gamma''''(w)) = 0$ and then entries
over the principal diagonal of
$\alpha_\gamma^{\ast}(P^{-1}\,\cdot\,dP)$ be zero. But this
determinant is in fact the differentiation of
$\det(\gamma'(s),\gamma''(s),\gamma'''(s))$ and for being zero it
is sufficient that we assume
$\det(\gamma'(s),\gamma''(s),\gamma'''(s))=1$. On the other hand,
we have
\begin{eqnarray*}
c' &=&(\gamma\circ\sigma)' = \sigma'\cdot(\gamma'\circ\sigma) \\
c'' &= & (\sigma')^2\cdot(\gamma''\circ\sigma)+\sigma''\cdot(\gamma'\circ\sigma)\\
c''' &=& (\sigma')^3\cdot(\gamma'''\circ\sigma)+ 3\,
\sigma'\sigma''\cdot(\gamma''\circ\sigma)+
\sigma'''\cdot(\gamma'\circ\sigma).
\end{eqnarray*}
Thus we conclude that
\begin{eqnarray*}
\det(c',c'',c''') &=&
\det(\sigma'\cdot(\gamma'\circ\sigma)\,,\,(\sigma')^2\cdot(\gamma''\circ\sigma)+\sigma''\cdot(\gamma'\circ\sigma)\,,\,\\
&& (\sigma')^3\cdot(\gamma'''\circ\sigma)+ 3\,
*\sigma'\sigma''\cdot(\gamma''\circ\sigma)+ \sigma'''\cdot
(\gamma'\circ\sigma)) \\
&&
\det(\sigma'\cdot(\gamma'\circ\sigma)\,,\,(\sigma')^2\cdot(\gamma''\circ\sigma)\,,\,(\sigma')^3\cdot
(\gamma'''\circ\sigma))\\
& = &(\sigma')^6 \cdot \det(\gamma'\circ\sigma\,,\,\gamma''\circ\sigma\,\,\gamma'''\circ\sigma)\\
& = &(\sigma')^6\,,
\end{eqnarray*}
The last expression specifies $\sigma$, namely the {\it special
affine arc length}, is defined as follows
$$\sigma:=\int_a^t\;\Big[\det(c'(u),c''(u),c'''(u))\Big]^{1/6}
du.$$\par
So $\sigma$ is a natural parameter under the action of special
affine transformations, that is, when $c$ is parameterized by
$\sigma$ then for each special affine transformation $A$, $A\circ
c$ is also parameterized by the same $\sigma$. Furthermore, every
curve parameterized by $\sigma$ up to special affine
transformations is introduced with the following invariants
\begin{eqnarray}
\chi_1=\det(c'',c''',c''''),\;\;\;\;\;\;\;\;\;
\chi_2=\det(c',c''',c'''').\label{(2)}
\end{eqnarray}
We call $\chi_1$ and $\chi_2$ the first and the second {\it
special affine curvatures} resp. Thus finally we have
\begin{eqnarray*}
\alpha_\gamma^{\ast}(P^{-1}\,.\,d P)=\left(
\begin{array}{ccc}
0 & 0
& \chi_1\\
1 & 0
& - \chi_2\\
0 & 1 & 0
\end{array}\right)d\sigma.
\end{eqnarray*}
\paragraph{Theorem 3.2} {\it Every space
curve of class ${\cal C}^5$ satisfying in condition (\ref{(1)})
under the action of special (unimodular) affine transformations is
determined by its natural equations $\chi_1=\chi_1(\sigma)$ and
$\chi_2=\chi_2(\sigma)$ of the first and second special affine
curvatures (\ref{(2)}) as functions (invariants) of the special
affine arc length.}
\paragraph{Theorem 3.3}{\it
Two space curves $c, \bar{c}: [a,b]\rightarrow{\bf R}^3$ of class
${\cal C}^5$ which satisfy in condition (\ref{(1)}) are special
affine equivalent if and only if $\chi_1^c=\chi_1^{\bar{c}}$ and
$\chi_2^c=\chi_2^{\bar{c}}$.}

\medskip \noindent {\it Proof:} The first side of the theorem is trivial with
respect to above descriptions. For the other side, we assume that
$c, \bar{c}$ are curves of class ${\cal C}^5$ satisfying (resp.)
in
\begin{eqnarray}
\det(c',c'',c''')>0,\hspace{1cm}
\det(\bar{c}',\bar{c}'',\bar{c}''')>0,\label{(3)}
\end{eqnarray}
that is, they are not plane curves. Also we suppose that functions
$\chi_1$ and $\chi_2$ are the same for these curves.\par
By changing the parameter to the natural parameter $\sigma$
(mentioned above), we obtain new curves $\gamma$ and
$\bar{\gamma}$ resp., that determinant (\ref{(3)}) will be equal
to 1. We prove that $\gamma$ and $\bar{\gamma}$ are special affine
equivalent and so there is a special affine transformation $A$
sech that $\bar{\gamma}=A\circ\gamma$. Then we have
$\bar{c}=A\circ c$ and the proof is complete.\par
At first, we replace the curve $\gamma$ with
$\delta:=\tau(\gamma)$ properly, in which case that $\delta$
intersects $\bar{\gamma}$ where $\tau$ is a translation defined by
translating one point of $\gamma$ to one point of $\bar{\gamma}$.
We correspond $t_0\in [a,b]$ to the intersection of $\delta$ and
$\bar{\gamma}$, thus $\delta(t_0)=\bar{\gamma}(t_0)$. One can find
a unique element $B$ of the general linear group ${\rm GL}(3,{\bf
R})$ such that maps the frame
$\{\delta'(t_0),\delta''(t_0),\delta'''(t_0)\}$ to the frame
$\{\bar{\gamma}'(t_0),\bar{\gamma}''(t_0),\bar{\gamma}'''(t_0)\}$.
So we have $B\circ\delta'(t_0)=\bar{\gamma}'(t_0)$,
$B\circ\delta''(t_0)=\bar{\gamma}''(t_0)$, and
$B\circ\delta'''(t_0)=\bar{\gamma}'''(t_0)$. $B$ is also an
element of the special linear group ${\rm SL}(3,{\bf R})$; since
we have
\begin{eqnarray*}
\det(\gamma'(t_0),\gamma''(t_0),\gamma'''(t_0))&=&\det(\delta'(t_0),\delta''(t_0),\delta'''(t_0))\;\;\;\;\textrm{and}\\
\det(\delta'(t_0),\delta''(t_0),\delta'''(t_0))&=&\det(B \circ
(\bar{\gamma}'(t_0),\bar{\gamma}''(t_0),\bar{\gamma}'''(t_0) )),
\end{eqnarray*}
and thus $\det(B)=1$. If we prove that $\eta:=B\circ\delta$ is
equal to $\bar{\gamma}$ on $[a,b]$, then by choosing $A=\tau\circ
B$, there will remain nothing for proof.

For curves $\eta$ and $\bar{\gamma}$ we have (resp.)
\begin{eqnarray*}
(\eta',\eta'',\eta''')'&=&(\eta',\eta'',\eta''')\left(
\begin{array}{ccc}
0 & 0 & \chi_1^{\eta} \\
1 & 0 & -\chi_2^{\eta} \\
0 & 1 & 0
\end{array}\right)\;\;\;\;\textrm{and}\\
(\bar{\gamma}',\bar{\gamma}'',\bar{\gamma}''')' &=&
(\bar{\gamma}',\bar{\gamma}'',\bar{\gamma}''')\left(
\begin{array}{ccc}
0 & 0 & \chi_1^{\bar{\gamma}} \\
1 & 0 & -\chi_2^{\bar{\gamma}} \\
0 & 1 & 0
\end{array}\right).
\end{eqnarray*}
Since $\chi_1$ and $\chi_2$ remain unchanged under special affine
transformations, so we have
$\chi_1^{\eta}=\chi_1^{\gamma}=\chi_1^{\bar{\gamma}}$ and
$\chi_2^{\eta}=\chi_2^{\gamma}=\chi_2^{\bar{\gamma}}$, therefore,
we conclude that $\eta$ and $\bar{\gamma}$ are solutions of
ordinary differential equation $Y''''+\chi_2\,Y''-\chi_1\,Y'=0$
where $Y$ depends to $t$. Because of the same initial conditions
\begin{eqnarray*}
&&\eta(t_0)=B\circ\delta(t_0)=\bar{\gamma}(t_0),\;\;\;\;\;\;\;\;\;
\eta'(t_0)=B\circ\delta'(t_0)=\bar{\gamma}'(t_0),\\
&&\eta''(t_0)=B\circ\delta''(t_0)=\bar{\gamma}''(t_0),\;\;\;\;\;
\eta'''(t_ )=B\circ\delta'''(t_0)=\bar{\gamma}'''(t_0),
\end{eqnarray*}
and the generalization of the existence and uniqueness theorem of
solutions, we have $\eta=\bar{\gamma}$ in a neighborhood of $t_0$
that can be extended to the whole $[a,b]$. \hfill\ $\diamondsuit$
\paragraph{Corollary 3.4} {\it The number of invariants of special
affine transformation group acting on ${\bf R}^3$ is two which is
the same with another results provided by other methods such as \cite{Gu}.}\\

The generalization of the affine classification of curves in an
arbitrary finite dimensional space has been discussed in
\cite{Na-Mah} and a complete set of invariants with a necessary
and sufficient condition of them for classifying curves up to
affine transformations has been derived.
%%%%%%%%%%%%%%%%%%%%%%%%%%%%%%%%%%%%%%%%%%%%%%%%%%%%%%%%%%%%%%%%%
\section{Geometric interpretations applied to physics and computer vision}
In the present section, the geometric interpretation of the first
and the second special affine curvatures and their applications in
physics and computer vision is discussed. Since every curve
parameterized with special affine arc length $\sigma$ and with
constant first and second affine curvatures $\chi_1$ and $\chi_2$
fulfilled in relation $\alpha'_c(\sigma)=\alpha_c(\sigma).(b)$,
for some $b\in {\rm sl}(3,{\bf R})$ via the right action of the
Lie algebra. In fact, we assumed that the action of the Lie group
be the left action \cite{Ol}. Whereof, Maurer-Cartan matrix of
${\rm SL}(3,{\bf R})$ is a base for Lie algebra ${\rm sl}(3,{\bf
R})$ and one can write
$\alpha'_c(\sigma)=\alpha_c(\sigma). \left(\begin{array}{ccc} 0 &
0 & \chi_1 \\ 1 & 0 & -\chi_2 \\
0 & 1 & 0  \end{array}\right)$.
By solving this first order equation, we obtain
$\alpha_c(\sigma)=\exp\Big(\sigma. \left(\begin{array}{ccc} 0 & 0
& \chi_1 \\ 1 & 0 & -\chi_2 \\
0 & 1 & 0  \end{array}\right)\Big)$ that, for different values of
$\chi_1$  and $\chi_2$ it has a different forms which we divide
these forms in the following cases:
%********************************************
\subsection*{I. {\bf The case $\chi_1=\chi_2=0$}.}
In this case, the curve is in the form
$\alpha_c(\sigma)=\left(\begin{array}{ccc} 1& 0
& 0 \\\sigma & 1 & 0 \\
\frac{1}{2}\,\sigma^2 & \sigma & 1  \end{array}\right)$. It is
clear that the first column of this matrix is $c'(\sigma)$ and so
we have
$c(\sigma)=K+(\sigma\,,\,\frac{1}{2}\,\sigma^2\,,\,\frac{1}{6}\,\sigma^3)$
for some constant $K\in{\bf R}^3$, that its image is analogous to
the image of {\it twisted cubic} \cite{Gr}. Also the image is
similar to the {\it Neil} or {\it semi-cubical} parabola's graph
\cite{Ku}. The projection of this space curve in the direction of
$z-$axis is a parabola in plane. This space curve is the simplest
curve in ${\bf R}^3$ under special affine transformations. Its
figure is a translation of Figure \ref{fig1}-(a) by constant $K$.
\paragraph{Theorem 4.1} {\it Space curves with zero special affine
curvatures are in the form of twisted cubic probably with some
translations.}\par
\begin{figure}[ht]
\epsfysize=4cm \centerline{\epsffile[0 0 819 446]{fig1.eps}}
\caption{(a) $\chi_1=\chi_2=0$. (b) $\chi_1=0,\chi_2>0$. }
\label{fig1}
\end{figure}
%
%***************************
\subsection*{II. {\bf The case $\chi_1=0$ and $\chi_2>0$}.}
In this case, we have
\begin{eqnarray*}
\alpha_c(\sigma)=\left(\begin{array}{ccc}
1& 0 & 0 \\
\frac{1}{\sqrt{\chi_2}}\,\sin(\sqrt{\chi_2}\,\sigma) &
\cos(\sqrt{\chi_2}\,\sigma) &
-\sqrt{\chi_2}\sin(\sqrt{\chi_2}\,\sigma) \\
-\frac{1}{\chi_2}\,(\cos(\sqrt{\chi_2}\,\sigma)-1) &
\frac{1}{\sqrt{\chi_2}}\,\sin(\sqrt{\chi_2}\,\sigma) &
\cos(\sqrt{\chi_2}\,\sigma)
\end{array}\right).
\end{eqnarray*}
So we obtain
$c(\sigma)=K+\Big(\sigma\,,\,-\frac{1}{\chi_2}\,\cos(\sqrt{\chi_2}\,\sigma)\,,\,
-\frac{1}{\chi_2\sqrt{\chi_2}}\,\sin(\sqrt{\chi_2}\,\sigma)+\frac{\sigma}{\chi_2}\Big)$
for $K\in{\bf R}^3$. The image of this curve is a translation of
Figure \ref{fig1}-(b) by constant $K$. Its projection in the
direction of $z-$axis is similar to the graph of function
$\cos(\sigma)$.
%
%***************************
\subsection*{III. {\bf The case $\chi_1=0$ and $\chi_2<0$}.}
If we use $|\chi_2|=-\chi_2$ to be the absolute value of $\chi_2$,
in the same way as the previous cases, we find that
\begin{eqnarray*} \alpha_c(\sigma)\!=\!\left(\!\!\begin{array}{ccc}
1& 0 & 0 \\
\frac{1}{\sqrt{|\chi_2|}}\,\sinh(\sqrt{|\chi_2|}\,\sigma) \!&\!
\cosh(\sqrt{|\chi_2|}\,\sigma) \!&\!
\sqrt{|\chi_2|}\sinh(\sqrt{|\chi_2|}\,\sigma) \\
\frac{1}{|\chi_2|}\,\{\cosh(\sqrt{|\chi_2|}\,\sigma)-1\} \!&\!
\frac{1}{\sqrt{|\chi_2|}}\,\sinh(\sqrt{|\chi_2|}\,\sigma) \!&\!
\cosh(\sqrt{|\chi_2|}\,\sigma)
\end{array}\!\!\right)\!.
\end{eqnarray*}
Thus $c(\sigma)= K + \Big(\sigma\, , \, \frac{1}{|\chi_2|}\, \cosh
(\sqrt{|\chi_2|}\, \sigma)\, , \,
 \frac{1}{|\chi_2|\sqrt{|\chi_2|}} \, \sinh ( \sqrt{|\chi_2|}\,
\sigma) - \frac{\sigma}{|\chi_2|} \Big)$ where $K$ is an element
of ${\bf R}^3$. Its image is drown in Figure \ref{fig2}-(a)
probably after a translation. Its $z-$axis projection is similar
to the graph of the function $\cosh(\sigma)$.
\begin{figure}[ht]
\epsfysize=4cm \centerline{\epsffile[0 0 598 323]{fig2.eps}}
\caption{(a) $\chi_1=0,\chi_2<0$. (b) $\chi_1>0,\chi_2=0$.}
\label{fig2}
\end{figure}
%***************************
\subsection*{IV. {\bf The case $\chi_1>0$ and $\chi_2=0$}.}
Under these conditions, the $\alpha_c(\sigma)$ is
\begin{eqnarray*}\left(\!\!\!\begin{array}{ccc}
\frac{2}{3}M\!+\!\frac{1}{3}R&
-\frac{1}{3}\chi_1^{1/3}(\sqrt{3}N\!+\!M\!-\!R) &
\frac{1}{3}\chi_1^{2/3}(\sqrt{3}N\!-\!M\!+\!R) \\
\frac{1}{3\chi_1^{1/3}}(\sqrt{3}N\!-\!M\!+\!R) &
\frac{2}{3}M+\frac{1}{3}R &
-\frac{1}{3}\chi_1^{1/3}(\sqrt{3}N\!+\!M\!-\!R) \\
-\frac{1}{3\chi_1^{2/3}}(\sqrt{3}N\!+\!M\!-\!R) &
\frac{1}{3\chi_1^{1/3}}(\sqrt{3}N\!-\!M\!+\!R) &
\frac{2}{3}M+\frac{1}{3}R
\end{array}\!\!\!\right),
\end{eqnarray*}
where
\begin{eqnarray*}\begin{array}{lcl}
M=\exp(-\frac{1}{2}\chi_1^{1/3}\sigma)\,\cos(\frac{\sqrt{3}}{2}\chi_1^{1/3}\sigma),&& R=\exp(\chi_1^{1/3}\sigma)\\
N=\exp(-\frac{1}{2}\chi_1^{1/3}\sigma)\,\sin(\frac{\sqrt{3}}{2}\chi_1^{1/3}\sigma).&&
\end{array}\end{eqnarray*}
Therefore with above conditions, we can write
\begin{eqnarray*}
c(\sigma)=K\!+\!\Bigg(\frac{1}{3\chi_1^{1/3}}\,(\sqrt{3}N\!-\!M\!+\!R),
\frac{1}{3\chi_1^{2/3}}\,(-\sqrt{3}N\!-\!M\!+\!R),
\frac{1}{3\chi_1}(2M\!+\!R)\Bigg)
\end{eqnarray*}
for some $K\in{\bf R}^{3}$. Its figure is similar to Figure
\ref{fig2}-(b).
%
%***************************
\subsection*{V. {\bf The case $\chi_1<0$ and $\chi_2=0$}.}
Such as case III, with use of $|\chi_2|=-\chi_2$, the conditions
of this case lead to the form of $\alpha_c(\sigma)$:
\begin{eqnarray*}
\left(\!\begin{array}{ccc} \frac{2}{3}M\!+\!\frac{1}{3}R &
\frac{1}{3}|\chi_1|^{1/3}(-\!\sqrt{3}N\!+\!M\!-\!R) &
\frac{1}{3}|\chi_1|^{2/3}(\sqrt{3}N\!+\!M\!-\!R) \\
\frac{1}{3|\chi_1|^{1/3}}(\sqrt{3}N\!+\!M\!-\!R) \!&\!
\frac{2}{3}M\!+\!\frac{1}{3}R \!&\!
\frac{-1}{3}|\chi_1|^{1/3}(\!\sqrt{3}N\!-\!M\!+\!R) \\
\frac{-1}{3|\chi_1|^{2/3}}(\!\sqrt{3}N\!-\!M\!+\!R) \!&\!
\frac{1}{3|\chi_1|^{1/3}}(\sqrt{3}N\!+\!M\!-\!R) \!&\!
\frac{2}{3}M\!+\!\frac{1}{3}R
\end{array}\!\right),
\end{eqnarray*}
where
\begin{eqnarray*}\begin{array}{lcl}
M=\exp(\frac{1}{2}|\chi_1|^{1/3}\sigma)\,\cos(\frac{\sqrt{3}}{2}|\chi_1|^{1/3/}\sigma),&& R=\exp(-|\chi_1|^{1/3}\sigma)\\
N=-\exp(\frac{1}{2}|\chi_1|^{1/3}\sigma)\,\sin(\frac{\sqrt{3}}{2}|\chi_1|^{1/3/}\sigma).&&
\end{array}\end{eqnarray*}
And so we have the following curve
\begin{eqnarray*}
c(\sigma)\!=\!K\!+\!\Big(\frac{1}{3|\chi_1|^{1/3}}(\sqrt{3}N\!+\!M\!-\!R),\frac{-1}{3|\chi_1|^{2/3}}(\sqrt{3}N\!
-\!M\!+\!R),\frac{-1}{3|\chi_1|})2M\!+\!R)\!\Big)
\end{eqnarray*}
for some $K\in{\bf R}^{3}$. Its shape is similar to Figure
\ref{fig3}-(a).
\begin{figure}[ht]
\epsfysize=4cm \centerline{\epsffile[0 0 360 384]{fig3.Eps}}
\caption{$\chi_1<0,\chi_2=0$.} \label{fig3}
\end{figure}
%
%***************************
\subsection*{VI. {\bf The case $\chi_1,\chi_2\neq 0$}.}
In this case, relations are not as simple as previous cases.
$\alpha_c(\sigma)$ in this general case, is in the following form
$\alpha_c(\sigma)=\left(\begin{array}{ccc}
B_{11}  & B_{12} & B_{13}\\
B_{21}  & B_{22} & B_{23} \\
B_{31}  & B_{32} & B_{33}
\end{array}\right)$
where for $1\leq i,j \leq 3$ the entries $B_{ij}$ for brevity are
given in Appendix at the end of the paper.

Thus the relation $c'(\sigma) = (B11,B21,B31)$ signifies the curve
$c(\sigma)$ by integrating of the coefficients with respect to
$\sigma$. Therefore we obtain $c(\sigma) = K + (T1, T2, T3)$ when
$K$ is an element of ${\bf R}^3$ and $T_i$~s are in the forms of
variables which are indicated in Appendix.

For different values of constants $\chi_1$, $\chi_2 \neq 0$, there
exist various curves and $c(\sigma)$ is a translation, contraction
or extraction of a curve in the form of these cases. Thus we have
different figures that each of which is similar to one of the
shapes given in Figure \ref{fig4}, (a)--(d).
\begin{figure}[ht]
\epsfysize=7.2cm \centerline{\epsffile[0 0 698 848]{fig4.eps}}
\caption{(a) $\chi_1,\chi_2>0$. \,\, (b) $\chi_1,\chi_2<0$.
\,\,(c) $\chi_1<0<\chi_2$.\,\, (d) $\chi_1>0>\chi_2$.}
\label{fig4}
\end{figure}
\paragraph{Corollary 4.2}{ \it
In general, every solution of $\alpha_c:{\bf R}\rightarrow {\rm
SL}(3,{\bf R})$ is provided by multiplying a special linear
transformation and a translation from an acquired curve in above
cases. In fact, the geometrical sense of above curves can
be explained as follows:\\
Each curve has two branches. The values of the first and second
special affine curvatures determine ``rotation quantities" of the
branches that by ascending (descending resp.) the values, each
branch's bend will increase (decrease resp.). Accordingly, the
definitions of $\chi_1$ and $\chi_2$ have geometric
interpretations as the usual terminology of curvatures.}\\

In the case of constant $\chi_1$ and $\chi_2$, by using Theorem
3.3, we can classify space curves in these cases via special
affine transformations and as a result we have the following
theorem:
\paragraph{Theorem 4.3}{\it Each curve of class ${\cal C}^5$ in ${\bf R}^3$
satisfied in condition (\ref{(1)}) with constant affine curvatures
$\chi_1$ and $\chi_2$, up to special affine transformations, is
the trajectory of a one--parameter subgroup of special
(unimodular) affine transformations, that is, a curve of cases
I-VI.}\\

Finally we give two applications of the classification of space
curves by the action of special affine transformations and Theorem
4.3:
\paragraph{Corollary 4.4}{ \it
In the physical sense, we may assume that each space curve ${\bf
X}:[a,b]\rightarrow{\bf R}^3$ is the trajectory of a particle with
a specified mass $m$" in ${\bf R}^3$ and in the view of an
observer, that is influenced under the effect of a force ${\bf
F}$. By the action of special affine transformations, particle's
path has two conservation laws: $({\bf X}''\times{\bf
X}''')\cdot{\bf X}''''$ and $({\bf X}'\times{\bf X}''')\cdot{\bf
X}''''$, that are, the first and the second special affine
curvatures. Therefore, by multiplying constant $m^3$ to theses
invariants, we find conservation laws as $({\bf F}\times{\bf
F}')\cdot{\bf F}''$ and $({\bf P}\times{\bf F}')\cdot{\bf F}''$
where ${\bf P}=m\cdot{\bf v}$ is the momentum of the particle. If
these invariants of the trajectory are constant, then the shape of
the motion is similar to one of the six cases mentioned in the
above theorem.}\\

The derived invariants in above, may have important applications
in astronomy, fluid mechanics, quantum, general relativity and
etc. A reason for this importance is that in these areas we deal
with the motion of a space particle and it may be of our interest
to investigate for symmetry properties and invariants of the
particle under rigid motions.
\paragraph{Corollary 4.5}{ \it In computer vision and image
processing, we may suppose that each space curve is one of the
characteristic curves on a 3-dimensional object, that are feasible
minimum segment curves that completely signify the object in the
viewpoint of an observer. Also, if by an effect provided by
(orientation-preserving) rotations and translations in ${\bf R}^3$
we change the position of a picture without any change in
characteristic lines, then these cur~es will be equivalent under
special affine transformation. If a characteristic line has
constant affine curvatures $\chi_1$ and $\chi_2$, then it will be
similar to one of the cases of curves mentioned in Theorem 4.3.}\\

For instance, these image invariants provide the most prominent
application fields in 3D medical imagery, including MRI,
ultrasound and CT data, in object recognition, symmetry and
differential invariant signatures of 3D shapes \cite{MP,MWTN,Ol1}.
%
%%%%%%%%%%%%%%%%%%%%%%%%%%%%%%%%%%%%%%%%%%%%
\section*{Appendix}
In case {\bf VI} of section 4, entries $B_{ij}$ ($1\leq i,j \leq
3$) are
\begin{eqnarray*}
B_{11}&=&\frac{1}{\Delta}\Big\{ A \Big(
72\,c_1^{1/3}\,\chi_1\,\chi_2-8\,c_1^{1/3}\,c_2\,\chi_2 -
144\, c_2\chi_1 - 192\,\chi_2^3-1296\,\chi_1^2 \\
&&+8\,c_1^{2/3}\,\chi_2^2 \big) +
B\Big(13.85640646\,c_1^{2/3}\,\chi_2^2 + 12.470765820
\,c_1^{1/3}\,\chi_1\,\chi_2\\
&&+ 13.85640646\,c_1^{1-3}\,c_2\,\chi_2 \Big) + D\Big(
-8\,c_1^{2/3}\,\chi_2^2 - 72\,c_2\,\chi_1 - 648\,\chi_1^2\\
&&- 96\,\chi_2^3 - 72\,b_1^{1/3}\,\chi_1\,\chi_2 -
8\,c_1^{1/3}\,c_2\,\chi_2 \Big)\Big\},\\
B_{21}&=& -B_{32}=-\frac{1}{\chi_1}\,B_{13}=
\frac{1}{\Delta}\Big\{ A\Big(-18\,c_1^{2/3}\,\chi_1 -
2\,c_1^{2/3}\,c_2 + 24\,c_1^{1/3}\,\chi_2^2\Big)\hspace{1cm} \\
&&+ B\Big(41.56921940\,c_1^{1/3}\chi_2^2 +
31.176874540c_1^{2/3}\chi_1 + 3.464101616c_1^{2/3}c_2\Big) \\
&&+ D\Big(-24\,c_1^{1/3}\,\chi_2^2 + 18\,c^{2/3}\,\chi_1 +
2\,c_1^{2/3}\,c_2\Big)\Big\},\\
B_{31}&=&\,c_1^{2/3}B_{12}= \frac{1}{\Delta}\Big\{
A\Big(-\,c_1^{4/3} - 12\,c_1^{2/3}\,\chi_2\Big)
+B\Big(-1.732050808\,c_1^{4/3} \\
&& + 20.7846097\,c_1^{2/3}\,\chi_2\Big) + D\Big(\,c_1^{4/3}+12\,c_1^{2/3}\,\chi_2\Big)\Big\},\\
B_{22}&=&B_{33}= \frac{1}{\Delta}\Big\{ A\Big(1296\,\chi_1^2 +
144\,c_2\,\chi_1 + 4\,c_1^{2/3}\,\chi_2^2 +
36\,c_1^{1/3}\,\chi_1\,\chi_2 \\
&&+ 4\,c_1^{1/3}\,c_2\,\chi_2 + 192\,\chi_2^3 \Big)+  B\Big(
62.353831080\,c_9^{1/3}\,\chi_1\,\chi_2 \\
&&+ 69.28203232\,c_1^{2/3}\,c_2\,\chi_2 -
69.28203232\,c_1^{2/3}\,\chi_2^2 \Big)+ D\Big(72\,c_2\,\chi_1 \\
&&+ 96\,\chi_2^3 + 648\,\chi_1^2 - 36\,c_1^{1/3}\,\chi_1\,\chi_2
-M 4\,c_1^{1/3}\,c_2\,\chi_2 -
4\,c_1^{2/3}\,\chi_2^2 \Big)\Big\},\\
F_{23}&=& \frac{-1}{\Delta}\Big\{A\Big(108 c_1^{1/3}\chi_1^2+24
c_1^{1/3}\chi_2^3 -6 c_1^{2/3}\chi_1\chi_2-2c_1^{1/3}c_2\chi_2+ 12 c_1^{1/3}c_2\chi_1\Big)\\
&&+ B\Big(10.39230485 c_1^{2/3}\chi_1\chi_2+
3.464101616 c_1^{1/3}c_2\chi_2+187.0614873\\
&&\times c_1^{1/3}\chi_1^2+ 41.5692194 c_1^{1/3}\chi_2^3+
20.7846097 c_1^{1/3}c_2\chi_1\Big)+D\Big(\!6c_1^{2/3}\chi_1\chi_2\\
&&-24 c_1^{1/3}\chi_2^3+ 2c_1^{2/3}c_2\chi_2 - 12
c_1^{2/3}c_2\chi_1- 108 c_1^{1/3}\chi_1^2\Big)\Big\},
\end{eqnarray*}
Which in the above relations we assumed that
\begin{eqnarray*}
&&\hspace{-0.7cm} c_1  = 108\,\chi_1 + 12\big(\sqrt{12\,\chi_2^3 + 81\,\chi_1^2}\big),\\
&&\hspace{-0.7cm} c_2 = \sqrt{12\,\chi_2^3 + 81\,\chi_1^2},\\
&&\hspace{-0.7cm} A = \exp\big(-0.08333333333
\frac{c_1^{2/3}-12\chi_1}{c_1^{1/3}}\,\sigma\big)\;
\cos\big(0.1443375673 \frac{c_1^{2/3} + 12\chi_1}{c_1^{1/3}}\,\sigma\big),\\
&&\hspace{-0.7cm} B = \exp\big(-0.08333333333
\frac{c_1^{2/3}-12\chi_1}{c_1^{1/3}}\,\sigma\big)\;
\sin\big(0.1443375673 \frac{c_1^{2/3} +
12\chi_1}{c_1^{1/3}}\,\sigma\big),\\
&&\hspace{-0.7cm} D = \exp\big(0.1666666667 \frac{c_1^{2/3}-12\chi_1}{c_1^{1/3}}\,\sigma\big),\\
&&\hspace{-0.7cm} \Delta = \big(c_1^{2/3}\,\chi_2 - 9
c_1^{1/3}\,\chi_1 - c_1^{1/3}\,c_2 -
12\,\chi_2^2\big)\big(c_1^{2/3} + 12\,\chi_2\big).
\end{eqnarray*}
Also $T_i$~s are in the following forms
\begin{eqnarray*} T_1&=& \frac{72}{\Delta'}\Big\{A
\Big(5832\chi_1^3\chi_2+864\chi_1\chi_2^4-24c_1^{1/3}\chi_2^5-81c_1^{2/3}\chi_1^3
+48c_2\chi_2^4 \\
&&-24c_1^{2/3}\chi_1\chi_2^3+648c_2\chi_1^2\chi_2-30c_1^{1/3}c_2\chi_2^2-9c_1^{2/3}
c_2\chi_1^2-2c_1^{2/3}c_2\chi_2^3\\
&&-270c_1^{1/3}\chi_1^2\chi_2^2\Big)+B\Big(
-41.56921939c_1^{1/3}\chi_2^5+140.2961154c_1^{2/3}\chi_1^3\\
&&+41.56921939c_1^{2/3}\chi_1\chi_2^3-467.6537182c_1^{1/3}\chi_1^2\chi_2^2
+ 15.58845727c_1^{2/3}c_2\chi_1^2\\
&&-51.96152424c_1^{1/3}c_2\chi_2^2+3.46101616c_1^{2/3}c_2\chi_2^3
\Big)+D\Big(432\chi_1\chi_2^4+24c_2\chi_2^4\\
&&+81c_1^{2/3}\chi_1^3+2916\chi_1^3\chi_2+24c_1^{1/3}\chi_
2^5+24c_1^{2/3}\chi_1\chi_2^3+270c_1^{1/3}\chi_1^2\chi_2^2\\
&&+324c_2\chi_1^2\chi_2+2c_1^{2/3}c_2\chi_2^3+9\,c_1^{2/3}\,c_2\,\chi_1^2
+30c_1^{1/3}c_2\chi_1\chi_2^2\Big)\Big\},\\
T_2&=& \frac{-1}{\Delta'}\Big\{
A\Big(864c_1^{2/3}\chi_1^2\chi_2+216c_1^{2/3}\,c_2\chi_1\chi_2+162c_1^{4/3}\chi_1^2
+18c_1^{4/3}c_2\chi_1\Big)\\
&&-B\Big(3367.1067714c_1^{2/3}\chi_1^2\chi_2-374.1229746c_1^{2/3}c_2\chi_1\chi_2\\
&&+280.59223086c_1^{4/3}\chi_1^2+31.176914544c_1^{4/3}c_2\chi_1\Big)-D\Big(162c_1^{4/3}\chi_1^2\\
&&-18c_1^{4/3}c_2\chi_1-216c_1^{2/3}c_2\chi_1\chi_2\Big)\Big\},\\
T_3&=&
\frac{-144}{\Delta'}\Big\{A\Big(12c_1^{1/3}\chi_2^4-324c_2\chi_1^2
-432\chi_1\chi_2^3-24c_2\chi_2^3-2916\chi_1^3+9c_1^{2/3}\\
&&\cdot\chi_1\chi_2^2+c_1^{2/3}c_2\chi_2^2+162c_1^{1/3}\chi_1^2\chi_2+18c_1^{1/3}c_2\chi_1\chi_2\Big)
+B\Big(20.78460970c_1^{1/3}\\
&&\cdot\chi_2^4-15.58845727c_1^{2/3}\chi_1\chi_2^2-1.732050808c_1^{2/3}c_2\chi_2^2
+280.5922309c_1^{1/3}\\
&&\cdot\chi_1^2\chi_2+31.17691454c_1^{1/3}c_2\chi_1\chi_2\Big)-D\Big(12c_1^{1/3}\chi_2^4
-9c_1^{2/3}\chi_1\chi_2^2-c_1^{2/3}c_2\chi_2^2\\
&&-162c_1^{1/3}\chi_1^2\chi_2-162c_2\chi_1^2-216\chi_1\chi_2^3-12c_2\chi_2^3
-1458\chi_1^3-18c_1^{1/3}c_2\chi_1\chi_2\Big)\Big\}.
\end{eqnarray*}
In the latter relations, $A$, $B$, $D$, $c_1$ and $c_2$ are the
same with relations mentioned in above and
\begin{eqnarray*}
\Delta'&=&\big[9c_1^{1/3}\chi_1+c_1^{1/3}c_2- c_1^{2/3}\chi_2+
12\chi_2^2\big]\big[9c_1^{1/3}\chi_1+
c_1^{1/3}c_2+c_1^{2/3}\chi_1+
12\chi_2^2\big]\\
&&\times \big[9c_1^{1/3}\chi_1+ c_1^{1/3}c_2-12\chi_2^2\big].
\end{eqnarray*}
%%%%%%%%%%%%%%%%%%%%%%%%%%%%%%%%%%%%%%%%%%%%
%

\end{document}